\documentclass[11pt,a4paper]{article}
\usepackage{amsmath}
\usepackage{amsfonts}
\usepackage[english]{babel}
\usepackage{amssymb,latexsym,dcolumn}
\usepackage{epsfig}

\usepackage{subfigure}

\newcommand{\CC}{\Bbb C}

\newcommand{\TT}{\Bbb T}

\newcommand{\er}{\mathbb{R}}
\newcommand{\cee}{\mathbb{C}}
\newcommand{\enn}{\mathbb{N}}
\newcommand{\zet}{\mathbb{Z}}

\newcommand{\bewijs}{\textsc{proof}}
\newcommand{\bol}{\hfill\square\\}

\newcommand{\til}{\tilde}

\newcommand{\spann}{\textrm{span}}
\newcommand{\el}{{\cal L}}

\newtheorem{st}{Theorem}
\newtheorem{prop}[st]{Proposition}

\newtheorem{lemma}[st]{Lemma}

\newtheorem{vb}[st]{Example}

\newtheorem{opm}[st]{Remark}
\newtheorem{gev}[st]{Corollary}

\title{Orthogonal Laurent polynomials on the unit circle and snake-shaped matrix factorizations}

\author{Ruym\'{a}n Cruz-Barroso\footnotemark[1]\ ,\quad Steven Delvaux\footnotemark[2]}
\date{}

\begin{document}

\maketitle
\renewcommand{\thefootnote}{\fnsymbol{footnote}}
\footnotetext[1]{Department of Computer Science, K.U.Leuven,
Celestijnenlaan 200A, B-3001 Leuven, Belgium. email:
ruyman.cruzbarroso\symbol{'100}cs.kuleuven.be. The work of this
author is partially supported by the research project MTM
2005-08571 of the Spanish Government and by the Fund of Scientific
Research (FWO), project ``RAM: Rational modeling: optimal
conditioning and stable algorithms'', grant $\sharp$G.0423.05 and
the Belgian Network DYSCO (Dynamical Systems, Control, and
Optimization), funded by the Interuniversity Attration Poles
Programme, initiated by the Belgian State, Science Policy Office.
The scientific responsibility rests with the author.}
\footnotetext[2]{Department of Mathematics, Katholieke
Universiteit Leuven, Celestijnenlaan 200B, B-3001 Leuven, Belgium.
Corresponding author, email:
steven.delvaux\symbol{'100}wis.kuleuven.be. The work of this
author is supported by the Onderzoeksfonds K.U.Leuven/Research
Fund K.U.Leuven.}


\begin{abstract}
Let there be given a probability measure $\mu$ on the unit circle
$\TT$ of the complex plane and consider the inner product induced
by $\mu$. In this paper we consider the problem of orthogonalizing
a sequence of monomials $\{z^{r_k}\}_k$, for a certain order of
the $r_k\in\zet$, by means of the Gram-Schmidt orthogonalization
process. This leads to a sequence of orthonormal Laurent
polynomials $\{\psi_k\}_k$. We show that the matrix representation
with respect to $\{\psi_k\}_k$ of the operator of multiplication
by $z$ is an infinite unitary or isometric matrix allowing a \lq
snake-shaped\rq\ matrix factorization. Here the \lq snake
shape\rq\ of the factorization is to be understood in terms of its
graphical representation via sequences of little line segments,
following an earlier work of S.~Delvaux and M.~Van Barel. We show
that the shape of the snake is determined by the order in which
the monomials $\{z^{r_k}\}_k$ are orthogonalized, while the \lq
segments\rq\ of the snake are canonically determined in terms of
the Schur parameters for $\mu$. Isometric Hessenberg matrices and
unitary five-diagonal matrices (CMV matrices) follow as a special
case of the presented formalism.


\end{abstract}

\section{Introduction}
\label{sectionintroduction}

\subsection{Isometric Hessenberg and unitary five-diagonal matrices}

In recent years, there has been a lot of research activity on the
topic of \emph{unitary five-diagonal matrices}, also known as
\emph{CMV matrices}. These matrices have been used by researchers
in various contexts, see e.g.\ \cite{CCG}, \cite{Ca}-\cite{CMV3},
\cite{Gol1}, \cite{Kil}, \cite{Nen}, \cite{BS}-\cite{BS2} and
\cite{Ve}.

Explicitly, the CMV matrix looks like
\begin{equation}\label{fullCMV}
{\cal C} = \left[\begin{array}{cccccccc}
\overline{\alpha_0} & \rho_0\overline{\alpha_1} & \rho_0\rho_1 & 0 & 0 & 0 & 0 & \ldots \\
\rho_0 & -\alpha_0\overline{\alpha_1} & -\alpha_0\rho_1 & 0 & 0 & 0 & 0 & \ldots \\
0 & \rho_1\overline{\alpha_2} & -\alpha_1\overline{\alpha_2} & \rho_2\overline{\alpha_3} & \rho_2\rho_3 & 0 & 0 & \ldots \\
0 & \rho_1\rho_2 & -\alpha_1\rho_2 & -\alpha_2\overline{\alpha_3} & -\alpha_2\rho_3 & 0 & 0 & \ldots \\
0 & 0 & 0 & \rho_3\overline{\alpha_4} & -\alpha_3\overline{\alpha_4} & \rho_4\overline{\alpha_5} & \rho_4\rho_5 & \ldots \\
0 & 0 & 0 & \rho_3\rho_4 & -\alpha_3\rho_4 & -\alpha_4\overline{\alpha_5} & -\alpha_4\rho_5 & \ldots \\
0 & 0 & 0 & 0 & 0 & \rho_5\overline{\alpha_6} & -\alpha_5\overline{\alpha_6} & \ldots \\
\vdots & \vdots & \vdots & \vdots & \vdots & \vdots & \vdots &
\ddots
\end{array}\right],
\end{equation}
where $\alpha_k$, $k=0,1,2,\ldots$ are complex numbers satisfying
$|\alpha_k|<1$ (the so-called \emph{Schur parameters} or
\emph{Verblunsky coefficients})
and $\rho_k:=\sqrt{1-|\alpha_k|^2}
\in (0,1]$ are the so-called \emph{complementary Schur
parameters}. The matrix ${\cal C}=(c_{i,j})_{i,j \geq
0}$\footnote{In the rest of the paper and for convenience with the
notation, we will label the rows and columns of any matrix
starting with index 0. As an example, the element $c_{1,1}$ in the
matrix (\ref{fullCMV}) will take the value
$-\alpha_0\overline{\alpha_1}$.} in \eqref{fullCMV} can be seen to
be unitary and five-diagonal, in the sense that $c_{i,j} = 0$
whenever $|i-j|>2$. More precisely, the nonzero entries of ${\cal
C}$ follow a kind of zigzag shape around the main diagonal.

The terminology \lq CMV matrix\rq\ for the matrix in \eqref{fullCMV}
originates from the book of Simon \cite{BS}, who named these
matrices after a 2003 paper by Cantero, Moral and Vel\'azquez
\cite{Ca}. But this terminology is far from historically correct,
since the latter paper \cite{Ca} is in fact a rediscovery of facts
which were already known by the numerical analysis community in the
early 1990's; a survey of these early results can be found in the
review paper by Watkins \cite{Wat}; see also \cite{BS2}.

In the present paper, we prefer to avoid such historical
discussions and we will therefore use the neutral term \lq unitary
five-diagonal matrix\rq\ to refer to these CMV matrices.

Unitary five-diagonal matrices have a number of interesting
features, including the statement proven in the literature that
(see further in this paper for more details) from all non-trivial
classes of unitary matrices, unitary five-diagonals have the
smallest bandwidth. Here the word \lq non-trivial\rq\ refers to
matrices which are not expressible as a direct sum of smaller
matrices.

While this statement about the minimal bandwidth is certainly
correct, it is a curious fact that this does \emph{not} imply that
unitary five-diagonal matrices are also \emph{numerically}
superior with respect to other non-trivial classes of
unitary/isometric matrices. For example, another class of
matrices which is often used in the literature is the class of
\emph{isometric Hessenberg matrices}, given explicitly by
\begin{equation}\label{fullHess}
{\cal H} = \left[\begin{array}{cccccc} \overline{\alpha_0} & \rho_0\overline{\alpha_1} & \rho_0\rho_1\overline{\alpha_2} & \rho_0\rho_1\rho_2\overline{\alpha_3} & \rho_0\rho_1\rho_2\rho_3\overline{\alpha_4} & \ldots  \\
\rho_0 & -\alpha_0\overline{\alpha_1} & -\alpha_0\rho_1\overline{\alpha_2} & -\alpha_0\rho_1\rho_2\overline{\alpha_3} & -\alpha_0\rho_1\rho_2\rho_3\overline{\alpha_4} & \ldots  \\
0 & \rho_1 & -\alpha_1\overline{\alpha_2} & -\alpha_1\rho_2\overline{\alpha_3} & -\alpha_1\rho_2\rho_3\overline{\alpha_4} & \ldots \\
0 & 0 & \rho_2 & -\alpha_2\overline{\alpha_3} & -\alpha_2\rho_3\overline{\alpha_4} & \ldots \\
0 & 0 & 0 & \rho_3 & -\alpha_3\overline{\alpha_4} & \ldots \\
 \vdots & \vdots & \vdots & \vdots & \vdots & \ddots
\end{array}\right].
\end{equation}
Note that the matrix in \eqref{fullHess} is of infinite dimension.
This matrix is called \emph{isometric} since its columns are
orthonormal; a similar property for the rows is not guaranteed.


In \eqref{fullHess} we use again the notation $\alpha_k$, $\rho_k$
to denote the Schur parameters and complementary Schur parameters,
respectively. These are the same numbers as in the matrix
\eqref{fullCMV}; see
further.

The matrix ${\cal H}=(h_{i,j})_{i,j \geq 0}$ in \eqref{fullHess} is
called \emph{Hessenberg} since $h_{i,j} = 0$ whenever $i-j\geq 2$.
Note however that the upper triangular part of this matrix is in
general dense.

Now the point is that isometric Hessenberg matrices as in
\eqref{fullHess} are known to be \emph{just as efficient} to
manipulate as unitary five-diagonal matrices!
Although this fact is known by numerical specialists, it seems
that it is not so well-known in part of the theoretical community.
Therefore, let us describe this now in somewhat more detail.

The naive idea would be that isometric Hessenberg matrices are \lq
inefficient\rq\ to work with since these matrices have a \lq
full\rq\ upper triangular part, in contrast to unitary
five-diagonal matrices. But this would be a too quick conclusion.
Having a better look at the problem, one can note that the upper
triangular part of an isometric Hessenberg matrix is \emph{rank
structured} in the sense that each submatrix that can be taken out
of the upper triangular part of such a matrix, has rank at most
equal to 1. This can be easily verified using e.g.\ the explicit
expressions of the entries of the matrix ${\cal H}$ in
\eqref{fullHess}.

Going one step further, one can note that the rank structure in the
upper triangular part of ${\cal H}$ is in fact a consequence of an
even more structural theorem. Denote with $G_{k,k+1}$ a \emph{Givens
transformation} (also called Jacobi transformation)
\begin{equation}\label{defgivens}
G_{k,k+1} = \left[\begin{array}{ccc} I_{k} & 0 & 0
\\ 0 & \tilde{G}_{k,k+1} & 0 \\ 0 & 0 & I \end{array}\right],
\end{equation}
where $I_{k}$ and $I$ denote identity matrices of sizes $k$ and
$\infty$, respectively, and where $\tilde{G}_{k,k+1}$ is a
2$\times$2 unitary matrix positioned in rows and columns
$\{k,k+1\}$. Thus the matrix $G_{k,k+1}$ differs from the identity
matrix only by its entries in rows and columns $\{k,k+1\}$. Givens
transformations can be considered as the most elementary type of
unitary matrices. They can be used as building blocks to construct
more general unitary matrices. Of interest for the present
discussion is the fact that any (infinite) isometric Hessenberg
matrix ${\cal H}$ allows a factorization as a product of Givens
transformations in the form
\begin{equation}\label{AGRfactorization} {\cal H} =
G_{0,1}G_{1,2}G_{2,3}G_{3,4}\ldots.\end{equation} This
factorization must be understood in the sense that the principal
$n\times n$ submatrices of ${\cal H}$ and $G_{0,1}G_{1,2}\ldots
G_{n-1,n}$ coincide for each $n$. This can be shown using only
some basic linear algebra \cite{GGMS,GVL}.

Applying this factorization to the matrix ${\cal H}$ in
\eqref{fullHess}, one can actually specify this result by noting
that the $k$th Givens transformation $G_{k,k+1}$ in
\eqref{AGRfactorization} must have nontrivial part given by
\begin{equation}\label{schurbuildingblock}\til{G}_{k,k+1} = \left[\begin{array}{cc} \overline{\alpha_k} & \rho_k \\
\rho_k & -\alpha_k\end{array}\right].\end{equation} In other
words, the \lq cosines\rq\ and \lq sines\rq\ of the Givens
transformations in \eqref{AGRfactorization} are nothing but the
Schur parameters and complementary Schur parameters, respectively.
This result was first established in the present context by Ammar,
Gragg and Reichel \cite{AGR}.

Incidently, note that the Givens transformations in
\eqref{schurbuildingblock} are of a special form in the sense that
they have real positive off-diagonal elements and determinant
$-1$.

We also note the following finite dimensional equivalent of
\eqref{AGRfactorization}: any unitary Hessenberg matrix ${\cal H}$
of size $n\times n$ allows a factorization in the form
\begin{equation}\label{AGRfactorizationbis}
{\cal H} = G_{0,1}G_{1,2}G_{2,3}G_{3,4}\ldots
G_{n-2,n-1},\end{equation} for suitable Givens transformations
$G_{k,k+1}$, $k = 0,1,\ldots,n-2$.

The main point of \eqref{AGRfactorizationbis} is that it shows
that unitary Hessenberg matrices of size $n$ can be compactly
represented using only $O(n)$ parameters, just as is the case for
unitary five-diagonal ones. Working with such an $O(n)$ matrix
representation, the eigenvalue problem for unitary Hessenberg
matrices can be solved numerically in a fast and accurate way; see
the end of Section \ref{sectionquadrature} for some references to
eigenvalue computation algorithms in the literature.
These algorithms can be
canonically expressed in terms of the matrix factorization
\eqref{AGRfactorizationbis}, i.e.,\ in terms of the Schur
parameters of the problem.

\subsection{Graphical representation}

In \cite{q351}, a graphical notation was introduced where matrix
factorizations with Givens transformations are represented via
sequences of little line segments. 

The graphical representation is obtained as follows. Let $A$ be
some arbitrary matrix (which will play no role in what follows)
and suppose that we update $A\mapsto G_{k,k+1}A$. This means that
the $k$th and $(k+1)$th row of $A$ are replaced by linear
combinations thereof, while the other rows of $A$ are left
unaltered. We can visualize this operation by drawing a vertical
line segment on the left of the two modified rows of $A$.

One can then apply this idea in an iterative way. For example,
when updating $A$ by means of an operation $A\mapsto
G_{k+1,k+2}G_{k,k+1}A$, one places first a vertical line segment
on the left of rows $k,k+1$ (this deals with the update $A\mapsto
G_{k,k+1}A$), and subsequently places a second vertical line
segment on the left of the former one, this time at the height of
rows $k+1,k+2$. We obtain in this way two successive vertical line
segments. Clearly, any number of Givens transformations can be
represented in such a way.

\begin{figure}[htbp]
\begin{center}
        \subfigure{\label{penta333brol}}\includegraphics[scale=0.3]{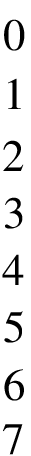}\hspace{10mm}
         \subfigure{\label{penta33brol}}\includegraphics[scale=0.3]{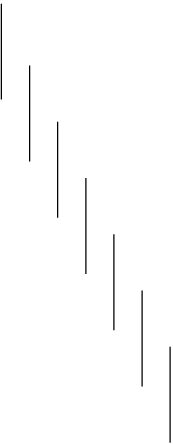}
\end{center}
\caption{The figure shows in a graphical way the decomposition as
a product of Givens transformations of the unitary Hessenberg
matrix ${\cal H}$ in \eqref{AGRfactorizationbis} with $n=8$.}
\label{fighess}
\end{figure}

Now the key point is that we \emph{identify} each $G_{k,k+1}$ with
its corresponding vertical line segment. We hereby make abstraction
of the matrix $A$ on whose rows these operations were assumed to
act. For example, the graphical representation of the factorization
\eqref{AGRfactorizationbis} with $n=8$ is shown in Figure
\ref{fighess}.

Concerning Figure \ref{fighess}, note that the top leftmost line
segment in this figure (which is assumed to be placed at \lq
height\rq\ 0 and 1; cf.\ the indices on the left of the figure)
corresponds to the leftmost factor $G_{0,1}$ in
\eqref{AGRfactorizationbis}. Similarly, the second line segment
corresponds to the factor $G_{1,2}$ in
\eqref{AGRfactorizationbis}, and so on. We again emphasize that
the line segments in Figure \ref{fighess} should be imagined as
\lq acting\rq\ on the rows of some (invisible) matrix $A$. See
\cite{q351,q352} for more applications of this graphical notation.

\begin{figure}[htbp]
\begin{center}
     \subfigure[]{\label{figpenta1}}\includegraphics[scale=0.3]{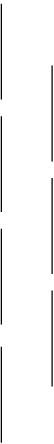}\hspace{20mm}
     \subfigure[]{\label{figpenta2}}\includegraphics[scale=0.3]{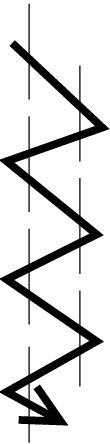}
\end{center}
\caption{The figure shows in a graphical way (a) the decomposition
as a product of Givens transformations of the unitary five-diagonal
matrix \eqref{givenspenta}, (b) the \lq snake shape\rq\ underlying
this decomposition.} \label{figpenta}
\end{figure}

It is known that also unitary five-diagonal matrices allow a
factorization as a product of Givens transformations. More precisely
\cite{AGR,BGE,Ca,Wat}, the matrix ${\cal C}$ in \eqref{fullCMV}
allows the factorization
$${\cal C} = \left[\begin{array}{cccccc} \overline{\alpha_0} & \rho_0 &  0 & 0 & 0 & \ldots \\
\rho_0 & -\alpha_0 & 0 & 0 & 0 & \ldots \\
0 & 0 & \overline{\alpha_2} & \rho_2 & 0  & \ldots \\
0 & 0 & \rho_2 & -\alpha_2 & 0  & \ldots \\
\vdots & \vdots  & \vdots & \vdots & \vdots & \ddots
\end{array}\right]\cdot
\left[\begin{array}{cccccc}
1 & 0 & 0 & 0 & 0 & \ldots \\
0 & \overline{\alpha_1} & \rho_1 & 0 & 0 &  \ldots \\
0 & \rho_1 & -\alpha_1 & 0 & 0 & \ldots \\
0 & 0 & 0 & \overline{\alpha_3} & \rho_3 & \ldots \\
0 & 0 & 0 & \rho_3 & -\alpha_3 &  \ldots \\
\vdots & \vdots & \vdots & \vdots & \vdots & \ddots
\end{array}\right],$$
which can be rewritten as
\begin{equation}\label{givenspenta} {\cal C} =
(\ldots G_{6,7}G_{4,5}G_{2,3}G_{0,1}) \cdot (G_{1,2}G_{3,4}G_{5,6}
\ldots),
\end{equation}
where the $G_{k,k+1}$ are again defined by \eqref{defgivens} and
\eqref{schurbuildingblock}. Again, this factorization must be
understood in the sense that the principal $n\times n$ submatrices
of ${\cal C}$ and $G_{n-2,n-1}\ldots G_{0,1} \cdot G_{1,2}\ldots
G_{n-1,n}$ for $n$ even or $G_{n-1,n}\ldots G_{0,1} \cdot
G_{1,2}\ldots G_{n-2,n-1}$ for $n$ odd coincide for each $n$. The
factorization \eqref{givenspenta} is represented graphically for
$n=8$ in Figure \ref{figpenta}.

Let us comment on Figure \ref{figpenta}. The leftmost series of line
segments in Figure \ref{figpenta1} corresponds to the leftmost
factor in \eqref{givenspenta}. The order in which these Givens
transformations are multiplied is clearly irrelevant; therefore we
are allowed to place them all graphically aligned along the same
vertical line. Similarly, the rightmost series of line segments in
Figure \ref{figpenta1} corresponds to the rightmost factor in
\eqref{givenspenta}. To explain Figure \ref{figpenta2}, imagine that
one moves from the top to the bottom of the graphical
representation. Then one can imagine a certain zigzag \lq snake
shape\rq\ underlying the factorization, which is shown in Figure
\ref{figpenta2}.

Note that in the above discussions we did not describe the way how
isometric Hessenberg and unitary five-diagonal matrices arise in
practice as matrix representations of a certain operator. At
present, it will suffice to know that they are matrix
representations of the \emph{operator of multiplication by $z$},
acting on a function space generated by a certain sequence of
orthonormal Laurent polynomials $\{\psi_k(z)\}_k$. This
orthonormal sequence
is obtained by applying the Gram-Schmidt
orthogonalization process to the sequence of monomials
\begin{equation}\label{vbhesspentaorder}
1,z,z^2,z^3,\ldots,\quad \textrm{and}\quad
1,z,z^{-1},z^2,z^{-2},\ldots,
\end{equation}
for the isometric Hessenberg and unitary five-diagonal case,
respectively. 

\subsection{Snake-shaped matrix factorizations}
\label{subsectionsnakeshapes}

The aim of this paper it to carry the above observations one step
further. We will show that with respect to a \emph{general}
sequence of orthonormal Laurent polynomials $\{\psi_k(z)\}_k$,
obtained by orthogonalizing a \emph{general} sequence of monomials
(satisfying some conditions to be described in detail in Section
\ref{subsectionlaurent}), the operator of multiplication by $z$ is
represented by an infinite unitary or isometric
matrix\footnote{Note of caution: we will also consider certain
cases where the subspace generated by the $\{\psi_k(z)\}_k$ is
\emph{not} invariant under the action of the operator of
multiplication by $z$. In such cases, the above statement has to
be formulated more carefully in order to make sure what the
meaning is of the matrix ${\cal S}$; actually this matrix needs
not be unitary nor isometric then. For a precise statement we
refer to the three cases distinguished at the beginning of Section
\ref{subsectionproof}, especially case 3.} allowing a snake-shaped
matrix factorization. We will use the latter term to denote an
infinite matrix product ${\cal S} = \prod_{k=0}^{\infty}
G_{k,k+1}$, where the factors under the $\prod$-symbol are
multiplied in a certain order. Here the \lq segments\rq\
$G_{k,k+1}$ of the snake are canonically fixed in terms of the
Schur parameters by means of \eqref{defgivens} and
\eqref{schurbuildingblock}, while the \lq shape\rq\ of the snake,
i.e., the order in which the $G_{k,k+1}$ are multiplied, will be
determined by the order in which the monomials have been
orthogonalized.

To fix the ideas, consider the sequence of monomials
\begin{equation}\label{vbgeneralorder}
1,z^{-1},z,z^{-2},z^2,z^3,z^{-3},z^{-4},z^4,z^5,\ldots.
\end{equation}  With respect to the
resulting sequence of orthonormal Laurent polynomials
$\{\psi_k(z)\}_k$ (see Section \ref{subsectionlaurent} for
details), the operator of multiplication by $z$ will be described
by a snake-shaped matrix factorization ${\cal S} = {\cal
S}^{(\infty)}$. We claim that this factorization is built by means
of the following recipe:
\begin{enumerate}
\item Considering the monomial $1=z^0$ in the position $0$ of
\eqref{vbgeneralorder}, we initialize ${\cal S}^{(0)}:= G_{0,1}$.
Then we apply the following procedure for $k\geq 1$:
\item If the $k$th monomial in \eqref{vbgeneralorder} has
a \emph{positive} exponent, we multiply the matrix with a new Givens
transformation on the \emph{right} by setting ${\cal S}^{(k)}:=
{\cal S}^{(k-1)}G_{k,k+1}$;
\item If the $k$th monomial in \eqref{vbgeneralorder} has a \emph{negative} exponent, we multiply the matrix with a new
Givens transformation on the \emph{left} by setting ${\cal
S}^{(k)}:= G_{k,k+1}{\cal S}^{(k-1)}$.
\end{enumerate} 

\begin{figure}[htbp]
\begin{center}
     \subfigure[]{\label{figsnake1}}\includegraphics[scale=0.3]{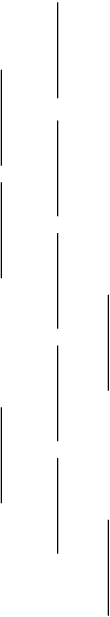}\hspace{20mm}
     \subfigure[]{\label{figsnake2}}\includegraphics[scale=0.3]{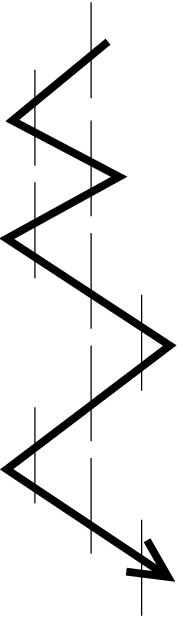}
\end{center}
\caption{The figure shows in a graphical way (a) the decomposition
as a product of Givens transformations of the matrix ${\cal S}$ in
\eqref{givenssnakeex}, (b) the \lq snake shape\rq\ underlying this
decomposition.} \label{figsnake}
\end{figure}

For the sequence of monomials \eqref{vbgeneralorder}, this recipe
gives rise to the following series of iterate matrices ${\cal
S}^{(k)}$:
$$\begin{array}{ll}
{\cal S}^{(0)} = G_{0,1}, &{\cal S}^{(1)} = G_{1,2} \cdot G_{0,1}, \\
{\cal S}^{(2)} = G_{1,2}\cdot G_{0,1}G_{2,3}, &{\cal S}^{(3)} =
G_{3,4}G_{1,2}\cdot
G_{0,1}G_{2,3} \\
{\cal S}^{(4)} = G_{3,4}G_{1,2}\cdot G_{0,1}G_{2,3}G_{4,5}, &{\cal
S}^{(5)} =
G_{3,4}G_{1,2}\cdot G_{0,1}G_{2,3}G_{4,5}G_{5,6}, \\
{\cal S}^{(6)} = G_{6,7}G_{3,4}G_{1,2}\cdot
G_{0,1}G_{2,3}G_{4,5}G_{5,6}, &\ldots.
\end{array}$$
This leads to the final matrix factorization
\begin{equation}\label{givenssnakeex}
{\cal S} = {\cal S}^{(\infty)} = (\ldots
G_{7,8}G_{6,7}G_{3,4}G_{1,2})\cdot(G_{0,1}G_{2,3}G_{4,5}G_{5,6}G_{8,9}G_{9,10}\ldots).
\end{equation}
This factorization 
is shown graphically in Figure \ref{figsnake}.

Let us comment on Figure \ref{figsnake}. The \lq snake\rq\ in this
figure was built by means of the following recipe:
\begin{enumerate}
\item Starting with a snake consisting of a single line segment $G_{0,1}$, we apply the following procedure for
$k\geq 1$:
\item If the $k$th monomial in \eqref{vbgeneralorder} has
a \emph{positive} exponent, the snake moves towards the
\emph{bottom right}, i.e., we add a new line segment on the bottom
right of the snake;
\item If the $k$th monomial in \eqref{vbgeneralorder} has a \emph{negative}
exponent, the snake moves towards the \emph{bottom left}, i.e., we
add a new line segment on the bottom left of the snake.
\end{enumerate}
Of course, this recipe is nothing but a direct translation of the
recipe that led us to the matrix factorization
\eqref{givenssnakeex}.

The reader should check that the above procedures are also valid
for the isometric Hessenberg and for the unitary five-diagonal
case (cf.~\eqref{vbhesspentaorder} and Figures \ref{fighess},
\ref{figpenta}).

\subsection{Outline and contributions of the paper}

The fact that the recipe in Section \ref{subsectionsnakeshapes}
leads to the correct matrix representation of the operator of
multiplication by $z$ with respect to the sequence of orthonormal
Laurent polynomials $\{\psi_k(z)\}_k$ will be shown in Section
\ref{sectionmainresult}. Our proof makes use of essentially three
facts: (i) an observation of Cruz-Barroso et al.~\cite{RO} (see
also Watkins \cite{Wat}) expressing the intimate connection
between orthonormal Laurent polynomials and Szeg\"o polynomials;
(ii) the well-known Szeg\"o recursion \cite{Sz}; and (iii) an
argument of Simon \cite{BS2} using \lq intermediary bases\rq\ in
the isometric Hessenberg case. The full proof is however rather
technical and requires some administrational book-keeping.

By factoring out a snake-shaped matrix product like
\eqref{givenssnakeex}, one can obtain explicit expressions for the
entries of the matrix, generalizing the expansions in
\eqref{fullCMV} and \eqref{fullHess}. This will be the topic of
Section \ref{sectionentrywiseexp}, where we will describe a
graphical rule for determining the zero pattern of the matrix ${\cal
S}$ as well as the shape of its non-zero elements.


Finally, in Section \ref{sectionquadrature} we will briefly
consider some connections between snake-shaped matrix
factorizations and Szeg\H{o} quadrature formulas. We will show
that the known results involving isometric Hessenberg and unitary
five-diagonal matrices can all be formulated in terms of a general
snake-shaped matrix factorization ${\cal S}$, extending an
observation of Ammar, Gragg and Reichel \cite{AGR}.

The remainder of this paper is organized as follows. Section
\ref{sectionmainresult} discusses some preliminaries about
sequences of orthogonal Laurent polynomials on the unit circle and
proves the main result about snake-shaped matrix factorizations.
Section \ref{sectionentrywiseexp} discusses the entry-wise
expansion of snake-shaped matrix factorizations. Finally, Section
\ref{sectionquadrature} considers the connection with Szeg\H{o}
quadrature formulas.

To end this introduction, let us discuss the main contributions of
this paper. It follows from the results presented here that
isometric Hessenberg and unitary five-diagonal matrices can be
considered as two extreme cases of a single mechanism, cf.\ the
discussion in Section \ref{subsectionsnakeshapes}. In this way we
obtain a unifying approach to some earlier results and estimates
in the literature, see e.g.\ \cite{CCG,CMV3,RO,BS}. In addition,
in the paper we provide graphical illustrations of the obtained
matrix factorizations. These graphics lead to additional insight,
explaining e.g.\ the term \lq snake-shaped matrix
factorization\rq. We feel that this might be an important
conceptual contribution in its own respect.

\section{Snake-shaped matrix factorizations: main result}
\label{sectionmainresult}

This section is devoted to the proof of our main result about
snake-shaped matrix factorizations, showing how these occur as the
matrix representation of the operator of multiplication by $z$
with respect to a sequence of orthonormal Laurent polynomials. We
start with some preliminaries.

\subsection{Sequences of orthogonal Laurent polynomials on the unit
circle} \label{subsectionlaurent}

In this first subsection we fix some notations and conventions
concerning orthogonal Laurent polynomials on the unit circle (see
\cite{CCG}, \cite{RP}-\cite{RO}). We denote by $\TT :=\{ z \in \CC
: |z|=1 \}$ the unit circle in the complex plane and by
$\Lambda:=\CC[z,z^{-1}]$ the complex vector space of Laurent
polynomials in the variable $z$. For a given order $n\in\enn$ and
an ordinary polynomial $p(z)=\sum_{k=0}^n c_{k}z^{k}$, we define
its dual as $p^*(z):=z^n\overline{p(1/\bar{z})}$, or explicitly
$p^*(z) = \sum_{k=0}^n \overline{c_{n-k}}z^{k}$. Here the bar
denotes complex conjugation.


Throughout the paper, we shall be dealing with a finite positive
non-discrete Borel measure $\mu$ supported on the unit circle
$\TT$ (which induces a measure on the interval $[-\pi,\pi]$ that
we also denote by $\mu$), normalized by the condition
$\int_{-\pi}^{\pi} d\mu(\theta)=1$, i.e., a probability measure.
As usual, the inner product induced by $\mu$ is given by
\begin{equation}\label{definnerproduct}\langle f,g
\rangle = \int_{-\pi}^{\pi} \overline{f\left( e^{i\theta} \right)}
g( e^{i\theta} ) d\mu(\theta),\end{equation} and the space of
quadratically integrable functions with respect to the inner
product \eqref{definnerproduct} is denoted as $L^{\mu}_2(\TT)$.

For our purposes, we start constructing a sequence of subspaces of
Laurent polynomials $\{ {\cal L}_n \}_{n=0}^{\infty}$ satisfying
$${\cal L}_0:=\textrm{span}\{ 1 \} \;\;\;,\;\; \textrm{dim}\left( {\cal L}_n \right) = n+1 \;\;\;,\;\; {\cal L}_n \subset {\cal L}_{n+1} \;\;,\;n \geq 1.$$
This can be done, by taking a sequence $\{ p_n \}_{n=0}^{\infty}$
of non-negative integers such that $p_0=0$, $0 \leq p_n \leq n$
and $s_n:=p_n-p_{n-1} \in \{0,1\}$ for all $n\geq 1$. In the
sequel, a sequence $\{ p_n \}_{n=0}^{\infty}$ satisfying these
requirements will be called a \emph{generating sequence}. Observe
that in this case both $\{ p_n \}_{n=0}^{\infty}$ and $\{ n-p_n
\}_{n=0}^{\infty}$ are non-negative non-decreasing sequences.
Then, set
$${\cal L}_n := \spann \left\{ z^j \;:\;-p_n \leq
j \leq n-p_n \right\}$$
and set ${\cal L}_{-1}:=\{ 0 \}$ to be the trivial subspace.
Observe that $\Lambda = \bigcup_{n=0}^{\infty} {\cal L}_n$ if and
only if $\displaystyle \lim_{n\rightarrow \infty} p_n =
\lim_{n\rightarrow \infty} (n-p_n) = \infty$ and that for all
$n\geq 1$,
$${\cal L}_{n} = \left\{ \begin{array}{lc}
{\cal L}_{n-1} \oplus \spann \{ z^{n-p_n} \} & \textrm{if} \;s_n=0, \\
{\cal L}_{n-1} \oplus \spann \{ z^{-p_n} \} & \textrm{if} \;s_n=1.
\end{array} \right.$$

Denote
\begin{equation}\label{defL}\el
:=\overline{\bigcup_{n=0}^{\infty} {\cal L}_n},\end{equation}
where $\overline{A}$ denotes the closure of $A$ with respect to
the norm induced by the inner product in $L^{\mu}_2(\TT)$. From
the fact that the Laurent polynomials form a dense subset in
$L^{\mu}_2(\TT)$, we have that ${\cal L} = L^{\mu}_2(\TT)$ if and
only if $\displaystyle \lim_{n\rightarrow \infty} p_n =
\lim_{n\rightarrow \infty} (n-p_n) = \infty$. If this condition is
violated, then ${\cal L}$ is a \emph{strict} subspace of
$L^{\mu}_2(\TT)$.

By applying the Gram-Schmidt orthogonalization procedure to ${\cal
L}_n$, an orthonormal basis $\{\psi_0(z), \ldots, \psi_n(z) \}$ can
be obtained. If we repeat the process for each $n \geq 0$, a
sequence $\{ \psi_n(z) \}_{n=0}^{\infty}$ of Laurent polynomials can
be obtained satisfying, for all $n,m\geq 0$:
\begin{enumerate}
\item $\psi_n(z) \in {\cal L}_n\setminus {\cal L}_{n-1},\ \ \quad \psi_0(z) \equiv 1$,
\item \vspace{-3mm}$\psi_n(z)$ has a real positive coefficient for the
power $\left\{
\begin{array}{ll}
z^{n-p_n} &\textrm{if } s_n=0 \\
z^{-p_n} &\textrm{if } s_n=1
\end{array}, \right.$
\item \vspace{-3mm}$\langle \psi_n(z),\psi_m(z) \rangle= \left\{
\begin{array}{ll}
0 &\textrm{if } n \neq m \\
1 &\textrm{if } n=m
\end{array}. \right.$
\end{enumerate}

This sequence will be called a {\em sequence of orthonormal
Laurent polynomials for the measure $\mu$ and the generating
sequence $\{ p_n \}_{n=0}^{\infty}$}.

Let us illustrate these ideas with three examples.

\begin{vb}
Consider the sequence of monomials given by (\ref{vbgeneralorder})
and the monomial $1=z^0$ in the position $0$. Then, the construction
of the sequence $\{ s_n \}_{n=1}^{\infty}$ is nothing but to take
$s_n=0$ if the $n$th monomial in (\ref{vbgeneralorder}) has a
positive exponent and $s_n=1$ if it is negative, whereas $p_n$
counts the number of negative monomials positioned up to $n$. Hence,
$\{ s_n \}_{n=1}^{\infty}=\{1,0,1,0,0,1,1,0,0,\ldots \}$ and $\{ p_n
\}_{n=0}^{\infty}=\{0,1,1,2,2,2,3,4,4,4,\ldots \}$.
\end{vb}

\begin{vb}
If $s_k = 0$ for all $k\geq 1$,
then ${\cal L}_n$ is the space of ordinary polynomials of degree
at most $n$. In this case the Gram-Schmidt orthogonalization
process is applied to the sequence of monomials $\{
1,z,z^2,z^3,\ldots\}$ and the resulting orthonormal Laurent
polynomials $\psi_n(z)$ are just the well-known orthonormal
\emph{Szeg\H{o} polynomials} $\varphi_n(z)$; see e.g.\ \cite{Sz}.
\end{vb}

\begin{vb}
If $s_k=k+1\mod 2$ for all $k\geq 1$,
then the Gram-Schmidt orthogonalization process is applied to the
sequence $\{1,z,z^{-1},z^2,z^{-2},\ldots\}$, where the monomials
$z^k$ and $z^{-k}$ occur in an alternating way. The resulting
sequence $\{ \psi_n(z) \}_{n=0}^{\infty}$ was firstly considered by
Thron in \cite{Th} and it is called the \emph{CMV basis} in
\cite{BS2}. The CMV basis can actually be expressed in terms of the
Szeg\H{o} polynomials as (see e.g. \cite{Ca,Cr,BS2,Th,Wat})
$$\varphi_0(z),\varphi_1(z),z^{-1}\varphi_2^*(z),z^{-1}\varphi_3(z),z^{-2}\varphi_4^*(z),z^{-2}\varphi_5(z),\ldots.$$
\end{vb}

In the general case, one has the following result.
\begin{lemma}\label{lema}
(Cruz-Barroso et al.~\cite{RO}; see also Watkins \cite{Wat}) The
family $\{ \psi_n(z) \}_{n=0}^{\infty}$ is the sequence of
orthonormal Laurent polynomials on the unit circle for a measure
$\mu$ and the ordering induced by the generating sequence $\{p_n
\}_{n=0}^{\infty}$, if and only if,
\begin{equation}\label{OLPversusSzego}
\psi_n(z) = \left\{\begin{array}{cc}
z^{-p_n}\varphi_n(z) & \textrm{if }s_n=0, \\
z^{-p_n}\varphi_n^*(z) & \textrm{if } s_n=1,
\end{array}\right.
\end{equation}
$\{\varphi_n(z)\}_{n=0}^{\infty}$ being the sequence of orthonormal
Szeg\H{o} polynomials for $\mu$.
\end{lemma}
$\bol$

Lemma \ref{lema} shows that the orthonormal Laurent polynomials
$\{\psi_n(z)\}_n$ are very closely related to the usual Szeg\"o
polynomials $\{\varphi_n(z)\}_{n}$ and their duals, and this for
\emph{any} choice of the generating sequence
$\{p_n\}_{n=0}^{\infty}$. We will need this result in what
follows.

\subsection{The main result}
\label{subsectionproof}

In this subsection we state and prove the main result of this
paper. Let $\{\psi_n(z)\}_{n=0}^{\infty}$ be the sequence of
orthonormal Laurent polynomials on the unit circle for the measure
$\mu$ and the ordering induced by the generating sequence
$\{p_n\}_{n=0}^{\infty}$. To distinguish them from the other
orthonormal sequences to be constructed in this section, we will
equip these Laurent polynomials with a superscript:
$\psi_n^{(0)}(z) := \psi_n(z)$. We will also find it convenient to
use the vectorial notation $\boldsymbol{\psi}^{(0)}(z) :=
(\psi^{(0)}_n(z))_{n=0}^{\infty}$. Thus, $\boldsymbol{\psi}^{(0)}$
is an infinite dimensional vector whose $n$th component is the
$n$th orthonormal Laurent polynomial $\psi^{(0)}_n$ ($n\geq 0$).

Let $M$ denote the operator of multiplication by $z$ on the space of
quadratically integrable functions with respect to the inner product
\eqref{definnerproduct}. Thus $M$ is defined by the action
$$ M:f(z)\mapsto zf(z) \;,\;\;f\in L^{\mu}_2(\TT).$$ Since we are working on the unit circle
$\TT$, the operator $M$ is actually \emph{unitary}.

Recall the notation ${\cal L}$ for the closure in $L^{\mu}_2(\TT)$
of the subspace generated by $\boldsymbol{\psi}^{(0)}(z)$. We
distinguish between three cases:

\begin{enumerate}
\item If $\lim_{n\rightarrow \infty}
p_n = \lim_{n\rightarrow \infty} (n-p_n) = \infty$, then $\el =
L^{\mu}_2(\TT)$. The sequence of orthonormal Laurent polynomials
$\boldsymbol{\psi}^{(0)}$ forms then a basis for $L^{\mu}_2(\TT)$
and the matrix representation of $M$ with respect to this basis is
an infinite unitary matrix ${\cal S}$, i.e., both the rows and
columns of this matrix are orthonormal.
\item If $\lim_{n\rightarrow \infty}
p_n < \infty$, then the sequence of orthonormal Laurent
polynomials $\boldsymbol{\psi}^{(0)}$ can be non-complete, but in
any way it will still generate a subspace of $L^{\mu}_2(\TT)$
which is invariant under the application of the operator $M$. We
can then define the operator $M\upharpoonright{\cal L}$, which is
the restriction of a unitary operator to an invariant subspace and
hence is isometric. The matrix representation of this operator
with respect to the basis $\boldsymbol{\psi}^{(0)}$ of $\el$ is
now an infinite isometric matrix ${\cal S}$, i.e., the columns of
this matrix are orthonormal. In fact, it is known that the
sequence $\boldsymbol{\psi}^{(0)}$ is complete in
$L^{\mu}_2(\TT)$, if and only if the so-called \emph{Szeg\"o
condition} fails, i.e., if $\sum_{j=0}^{\infty}
|\alpha_j|^2=\infty$. In that case the matrix ${\cal S}$ is
actually \emph{unitary} since $M\upharpoonright{\cal L}=M$.

\item If $\lim_{n\rightarrow \infty} (n-p_n) < \infty$\footnote{We
thank the referee for pointing our attention to this case, and for
providing us with the modifications that have to be made for it.},
then the sequence of orthonormal Laurent polynomials
$\boldsymbol{\psi}^{(0)}$ can be non-complete, and in that case it
generates a subspace of $L^{\mu}_2(\TT)$ which is \emph{not}
invariant under the application of the operator $M$. However, we
can now still consider the operator $PM\upharpoonright{\cal L}$
where $P$ is the orthogonal projection operator of
$L^{\mu}_2(\TT)$ onto $\el$. The matrix representation ${\cal S}$
of this operator with respect to the basis
$\boldsymbol{\psi}^{(0)}$ of $\el$ is now not necessarily unitary
neither isometric. Actually, it holds that the \emph{rows} of this
matrix are orthonormal. This follows by noticing that the
transpose of the matrix ${\cal S}$ occurs as a matrix
representation in the previous case and hence is isometric.
\end{enumerate}

Note that in each of the above three cases, the infinite matrix
${\cal S}$ has its entries given by
\begin{eqnarray}
\nonumber {\cal S} &=& [\langle \psi^{(0)}_i(z),z\psi^{(0)}_j(z)
\rangle]_{i,j=0}^{\infty} \\ \label{firstbasis} &=:& \langle
\boldsymbol{\psi}^{(0)}(z),z\boldsymbol{\psi}^{(0)}(z) \rangle,
\end{eqnarray}
where the inner product is defined in \eqref{definnerproduct}.
Here the expression on the second line should be regarded as a
compact vectorial notation of the line above.

Now we are in position to prove the following result. We will do
this by using a modification of an argument of Simon \cite[third
proof of Theorem 10.1]{BS2} for the isometric Hessenberg case. The
main ingredient of the proof will be the well-known
\emph{Szeg\H{o} recursion}, expressed in the form (see e.g.\
\cite{Sz})
\begin{equation}\label{Szegomatrixrecursion}
\left[\begin{array}{c} z\varphi_k(z) \\ \varphi_{k+1}^{*}(z)
\end{array}\right] = \left[\begin{array}{cc} \overline{\alpha_k} & \rho_k \\ \rho_k & -\alpha_k
\end{array}\right] \left[\begin{array}{c} \varphi_k^{*}(z) \\ \varphi_{k+1}(z)
\end{array}\right],
\end{equation}
where $\varphi_k(z)$ and $\varphi_{k}^{*}(z)$ denote the
orthonormal Szeg\H{o} polynomial of degree $k$ and its dual
respectively. Note that the coefficient matrix in
\eqref{Szegomatrixrecursion} is nothing but the nontrivial part
\eqref{schurbuildingblock} of the Givens transformation
$G_{k,k+1}$.

\begin{st}\label{maintheorem} Let $\{\psi_n(z)\}_{n=0}^{\infty}$ be
the sequence of orthonormal Laurent polynomials on the unit circle
for a measure $\mu$ and the ordering induced by the generating
sequence $\{p_n \}_{n=0}^{\infty}$. Then the matrix ${\cal S}$ in
\eqref{firstbasis} can be factored into a snake-shaped matrix
factorization ${\cal S}={\cal S}^{(\infty)}$, constructed by the
recipe given in Section \ref{subsectionsnakeshapes}. The
factorization must be understood in the sense that the principal
$n\times n$ submatrices of ${\cal S}^{(n-1)}$ and ${\cal S}$
coincide for all $n$.
\end{st}

\bewijs. We will construct a sequence of intermediary bases
$\boldsymbol{\psi}^{(k)}$ for the subspace
$\el=\overline{\spann\{\psi^{(0)}_j\}_{j=0}^{\infty}}$, $k\geq 1$,
in such a way that for each $k$, there exists an index
$l\in\{0,1,\ldots,k-1\}$ such that $\boldsymbol{\psi}^{(k)}$ is
the same as $\boldsymbol{\psi}^{(l)}$, except for a change in the
$(k-1)$th and $k$th components. These intermediary bases will
serve to factorize the matrix ${\cal S}$. For example, note that
\eqref{firstbasis} can be rewritten as
\begin{eqnarray}
\label{secondbasisaux} {\cal
S} & = &  \langle \boldsymbol{\psi}^{(0)}(z),z\boldsymbol{\psi}^{(0)}(z) \rangle \\
\label{secondbasis} & = & \langle
\boldsymbol{\psi}^{(0)}(z),\boldsymbol{\psi}^{(1)}(z) \rangle
\cdot\langle \boldsymbol{\psi}^{(1)}(z),z\boldsymbol{\psi}^{(0)}(z)
\rangle,
\end{eqnarray}
for any choice of the basis $\boldsymbol{\psi}^{(1)}$ of $\el$.
Indeed, the $j$th column of the matrix \eqref{secondbasis} is
obtained by expressing the orthogonal projection on $\el$ of the
function $z\psi^{(0)}_{j}(z)$ in terms of the basis
$\boldsymbol{\psi}^{(1)}(z)$, which is then in its turn expressed
in terms of the basis $\boldsymbol{\psi}^{(0)}(z)$. Obviously this
gives the same result as directly expressing the orthogonal
projection on $\el$ of $z\psi^{(0)}_{j}(z)$ in terms of the basis
$\boldsymbol{\psi}^{(0)}(z)$, i.e., it equals the $j$th column of
\eqref{secondbasisaux} (we recall again our convention with the
notation: $j\geq 0$).

Note that instead of \eqref{secondbasis} we could also have
written a slightly modified version of it:
\begin{eqnarray}
\nonumber {\cal
S} & = &  \langle \boldsymbol{\psi}^{(0)}(z),z\boldsymbol{\psi}^{(0)}(z) \rangle \\
\nonumber & = & \langle
\boldsymbol{\psi}^{(0)}(z),z\boldsymbol{\psi}^{(1)}(z) \rangle
\cdot\langle
z\boldsymbol{\psi}^{(1)}(z),z\boldsymbol{\psi}^{(0)}(z) \rangle \\
\label{thirdbasis} & = & \langle
\boldsymbol{\psi}^{(0)}(z),z\boldsymbol{\psi}^{(1)}(z) \rangle
\cdot\langle \boldsymbol{\psi}^{(1)}(z),\boldsymbol{\psi}^{(0)}(z)
\rangle
\end{eqnarray}
for any choice of the basis $\boldsymbol{\psi}^{(1)}$ of $\el$ and
where we have used the general fact that $\langle zf(z),zg(z)
\rangle = \langle f(z),g(z) \rangle$ for any functions
$f,g:\TT\to\cee$, which follows from \eqref{definnerproduct} and
the fact that $z\in\TT$. The choice between \eqref{secondbasis}
and \eqref{thirdbasis} will depend on the fact whether $s_1=0$ or
$s_1=1$, respectively; see further.

The point will now be to make a good choice for the intermediary
bases $\boldsymbol{\psi}^{(k)}$. For example, the \lq good\rq\
choice for $\boldsymbol{\psi}^{(1)}$ will be the one for which one
of the factors in \eqref{secondbasis} (or \eqref{thirdbasis}) equals
the Givens transformation $G_{0,1}$, while the other factor is of
the form
$$\left[\begin{array}{cc} 1 & 0 \\ 0 & * \end{array}\right],$$
where $*$ denotes an irrelevant submatrix (which is actually of
infinite dimension). Explicitly, the basis $\boldsymbol{\psi}^{(1)}$
is given by $\boldsymbol{\psi}^{(1)}_0=z^{1-2s_1}$,
$\boldsymbol{\psi}^{(1)}_1=z^{-s_1}[z^{s_1}
\boldsymbol{\psi}^{(0)}_{1}]^*$ and
$\boldsymbol{\psi}^{(1)}_k=\boldsymbol{\psi}^{(0)}_k$ for all $k\geq
2$. Repeating this idea inductively for all subsequent bases
$\boldsymbol{\psi}^{(k)}$ will ultimately lead to the decomposition
of ${\cal S}$ as an infinite product of Givens transformations.

Let us now formalize these ideas. We work with the induction
hypothesis that after the $k$th step, $k\geq 0$\footnote{This
procedure also works for $k=0$ provided that we set ${\cal
S}^{(-1)}:=I$ and $s_0=0$ or $s_0=1$, since either choice will give
the same result.}, we have decomposed the matrix ${\cal S}$ as
\begin{equation}\label{inductionSX}
\begin{array}{c}
{\cal S} = S^{(k-1)}X^{(k)},\quad \textrm{ if } s_k=0, \\
{\cal S} = X^{(k)}S^{(k-1)},\quad \textrm{ if } s_k=1,
\end{array}
\end{equation}
where ${\cal S}^{(k-1)}$ is the $(k-1)$th iterate matrix of the
snake-shaped matrix factorization ${\cal S}^{(\infty)}$ (cf.\ the
construction in Section \ref{subsectionsnakeshapes}), while
$X^{(k)}$ equals the identity matrix in its first $k$ rows and
columns, i.e.,
\begin{equation}\label{inductionX1} X^{(k)} = \left[\begin{array}{cc}
I_{k} & 0  \\
0 & *
\end{array}\right],
\end{equation}
with $I_{k}$ the identity matrix of size $k$. We also assume by
induction that
\begin{equation}\label{inductionX2}
X^{(k)} = \langle
\boldsymbol{\psi}^{(l)}(z),z\boldsymbol{\psi}^{(m)}(z) \rangle,
\end{equation}
where $l,m$ are certain indices in $\{0,1,\ldots,k\}$ with at
least one of them equal to $k$ (we could actually give explicit
expressions for $l,m$ but will not need these in what follows).
Note that by combining the hypotheses \eqref{inductionX1} and
\eqref{inductionX2}, we deduce that $\psi_i^{(l)}(z) =
z\psi_i^{(m)}(z)$ for all $i\in\{0,1,\ldots,k-1\}$. In addition,
we have the following induction hypothesis on the $k$th components
of $\boldsymbol{\psi}^{(l)}$ and $\boldsymbol{\psi}^{(m)}$:
\begin{equation}\label{inductionkthcomponent} \begin{array}{c}\psi^{(l)}_{k}(z) = z^{-p_k}\varphi_{k}^*(z),\\
\psi^{(m)}_{k}(z) = z^{-p_k}\varphi_{k}(z),
\end{array}
\end{equation}
where $\varphi_{k}(z)$ denotes the orthonormal Szeg\H{o} polynomial
of degree $k$.

Identities \eqref{inductionSX} and \eqref{inductionX1} imply the
coincidence of the principal $k\times k$ submatrices of
$S^{(k-1)}$ and $S$, as the theorem states. Thus, to prove the
theorem we simply must show that, given all the above induction
hypotheses, we can now come to the induction step $k\mapsto k+1$.
To this end, we should try to peel off a new Givens transformation
$G_{k,k+1}$ from the matrix ${\cal S}$. Assume that the first $k$
intermediary bases $\boldsymbol{\psi}^{(1)}, \ldots,
\boldsymbol{\psi}^{(k)}$ of $\el$ have already been constructed.
We want to define the next intermediary basis
$\boldsymbol{\psi}^{(k+1)}$. We distinguish between two cases:
\begin{enumerate}
\item If $s_{k+1}=0$, we define $\boldsymbol{\psi}^{(k+1)}$ to be the same as
$\boldsymbol{\psi}^{(l)}$, except for its $k$th and $(k+1)$th
components. More precisely, we set
\begin{eqnarray}
\label{updatecase11}\left[\begin{array}{c} \psi^{(k+1)}_{k}(z)\\
\psi^{(k+1)}_{k+1}(z)
\end{array}\right] & := & \tilde{G}_{k,k+1}\left[\begin{array}{c} \psi^{(l)}_{k}(z)\\ \psi^{(l)}_{k+1}(z)
\end{array}\right]\\
\label{updatecase12} & = & \tilde{G}_{k,k+1}\cdot z^{-p_{k+1}}\left[\begin{array}{c} \varphi_{k}^*(z)\\
\varphi_{k+1}(z)
\end{array}\right]  \\
\label{updatecase13} & = & z^{-p_{k+1}}\left[\begin{array}{c} z\varphi_{k}(z)\\
\varphi_{k+1}^*(z)\end{array}\right].
\end{eqnarray}
Here the second equality follows from the first lines of
\eqref{inductionkthcomponent} and \eqref{OLPversusSzego} (recall
that the $(k+1)$th component of $\boldsymbol{\psi}^{(l)}$ has not
been changed yet with respect to $\boldsymbol{\psi}^{(0)}$), and
from the fact that $p_{k+1} = p_k$ by assumption. On the other
hand, the third equality is nothing but the Szeg\H{o} recursion
\eqref{Szegomatrixrecursion}.

Then, we can factorize \eqref{inductionX2} as
\begin{eqnarray}
\nonumber X^{(k)} &=& \langle
\boldsymbol{\psi}^{(l)}(z),z\boldsymbol{\psi}^{(m)}(z) \rangle\\
\nonumber &=& \langle
\boldsymbol{\psi}^{(l)}(z),\boldsymbol{\psi}^{(k+1)}(z) \rangle\cdot
\langle
\boldsymbol{\psi}^{(k+1)}(z),z\boldsymbol{\psi}^{(m)}(z) \rangle\\
\nonumber &=& G_{k,k+1}\cdot\langle
\boldsymbol{\psi}^{(k+1)}(z),z\boldsymbol{\psi}^{(m)}(z) \rangle\\
\label{updatecase14} &=:& G_{k,k+1}X^{(k+1)},
\end{eqnarray}
where the third equality follows from \eqref{updatecase11}, and
where the matrix $X^{(k+1)}$ in the fourth equality now equals the
identity matrix in its first $k+1$ rows and columns. The latter
follows by the induction hypothesis for the first $k$ rows and
columns (rows and columns $0$ to $k-1$), and from the fact that,
by the first line of \eqref{updatecase13} and the second line of
\eqref{inductionkthcomponent}, we have
$$\psi^{(k+1)}_k(z) = z\cdot z^{-p_{k+1}}\varphi_{k}(z) = z\psi^{(m)}_k(z),$$ implying that also
the $k$th column of the matrix $X^{(k+1)}$ has all its entries
equal to zero, except for the diagonal entry which equals one.
From the fact that the $k$th row of the matrix $X^{(k+1)}$ is a
vector with norm at most one and with one of its entries equal to
one, it then follows that also the $k$th \emph{row} has all its
entries equal to zero, except for the diagonal entry.

We can then replace the index $l$ by its new value $k+1$. We have
already checked that the induction hypotheses \eqref{inductionX1}
and \eqref{inductionX2} are inherited in this way as $k\mapsto k+1$.
Also the hypothesis \eqref{inductionkthcomponent} can be easily
checked to remain valid in this way, by virtue of the second line of
\eqref{updatecase13} and the first line of \eqref{OLPversusSzego}.
Finally, we have to check that \eqref{inductionSX} remains also
valid. To prove this, we use \eqref{inductionSX},
\eqref{updatecase14}, the construction of ${\cal S}^{(k)}$ in
Section \ref{subsectionsnakeshapes} and we distinguish between two
cases:
\begin{enumerate}
\item If $s_k=0$ then
\begin{eqnarray*}
{\cal S}& = & {\cal S}^{(k-1)}X^{(k)} \\ &=& {\cal
S}^{(k-1)}G_{k,k+1}X^{(k+1)} \\
&=:& {\cal S}^{(k)}X^{(k+1)}. \end{eqnarray*}
\item If $s_k=1$ then
\begin{eqnarray*} {\cal
S} & = & X^{(k)}{\cal S}^{(k-1)}\\  &=& G_{k,k+1}X^{(k+1)}{\cal
S}^{(k-1)}\\ &=& G_{k,k+1}{\cal
S}^{(k-1)}X^{(k+1)}\\
&=:& {\cal S}^{(k)}X^{(k+1)},\end{eqnarray*} where we have used the
commutativity of ${\cal S}^{(k-1)}$ and $X^{(k+1)}$ since these
matrices have a complementary zero pattern (the former equals the
identity matrix except for its first $(k+1)\times (k+1)$ block,
while the latter is precisely the identity matrix there, cf.\
\eqref{inductionX1}).
\end{enumerate}

\item If $s_{k+1}=1$, we define $\boldsymbol{\psi}^{(k+1)}$ to be the same as
$\boldsymbol{\psi}^{(m)}$, except for its $k$th and $(k+1)$th
components. More precisely, we set
\begin{eqnarray}
\label{updatecase21}\left[\begin{array}{c} \psi^{(k+1)}_{k}(z)\\
\psi^{(k+1)}_{k+1}(z)
\end{array}\right] & := & \tilde{G}_{k,k+1}^{-1}\left[\begin{array}{c} \psi^{(m)}_{k}(z)\\ \psi^{(m)}_{k+1}(z)
\end{array}\right]\\
\label{updatecase22} & = &
\tilde{G}_{k,k+1}^{-1}\cdot z^{-p_{k+1}}\left[\begin{array}{c} z\varphi_{k}(z)\\
\varphi^*_{k+1}(z)
\end{array}\right] \\
\label{updatecase23} &=& z^{-p_{k+1}}\left[\begin{array}{c} \varphi_{k}^*(z)\\
\varphi_{k+1}(z)
\end{array}\right],
\end{eqnarray}
where we have used the second lines of
\eqref{inductionkthcomponent} and \eqref{OLPversusSzego} (recall
that the $(k+1)$th component of $\boldsymbol{\psi}^{(m)}$ has not
been changed yet with respect to $\boldsymbol{\psi}^{(0)}$), the
fact that $p_{k+1}=p_k+1$ by assumption and the Szeg\H{o}
recursion \eqref{Szegomatrixrecursion}.

Then, we can factorize \eqref{inductionX2} as
$$\begin{array}{ccl}
\nonumber X^{(k)} &=& \langle
\boldsymbol{\psi}^{(l)}(z),z\boldsymbol{\psi}^{(m)}(z) \rangle\\
\nonumber &=& \langle
\boldsymbol{\psi}^{(l)}(z),z\boldsymbol{\psi}^{(k+1)}(z)
\rangle\cdot \langle
z\boldsymbol{\psi}^{(k+1)}(z),z\boldsymbol{\psi}^{(m)}(z) \rangle\\
\nonumber &=& \langle
\boldsymbol{\psi}^{(l)}(z),z\boldsymbol{\psi}^{(k+1)}(z)
\rangle\cdot \langle
\boldsymbol{\psi}^{(k+1)}(z),\boldsymbol{\psi}^{(m)}(z) \rangle\\
\nonumber &=& \langle
\boldsymbol{\psi}^{(l)}(z),z\boldsymbol{\psi}^{(k+1)}(z)
\rangle\cdot G_{k,k+1}\\
&=:& X^{(k+1)}G_{k,k+1},
\end{array}$$
where the fourth step follows from \eqref{updatecase21}.

It is easy to check again that $X^{(k+1)}$ equals the identity
matrix in its first $k+1$ rows and columns by using the induction
hypothesis for the first $k$ rows and columns ($0,1,\ldots,k-1$)
and from the first lines of \eqref{updatecase23} and
\eqref{inductionkthcomponent} for the $k$th row and column.

We can then replace the index $m$ by its new value $k+1$. If follows
from the above discussion that the induction hypotheses
\eqref{inductionX1} and \eqref{inductionX2} are inherited in this
way as $k\mapsto k+1$. Also the hypothesis
\eqref{inductionkthcomponent} goes through, by virtue of the second
lines of \eqref{OLPversusSzego} and \eqref{updatecase23}. Finally,
the proof that also \eqref{inductionSX} goes through can be proven
in a completely similar way as in the previous case.
\end{enumerate}
We have now completely established the induction hypothesis
$k\mapsto k+1$, hereby ending the proof of Theorem
\ref{maintheorem}. $\bol$

\section{Entry-wise expansion of a snake-shaped matrix factorization}
\label{sectionentrywiseexp}

In this section we discuss the entry-wise expansion of a
snake-shaped matrix factorization ${\cal S}$, hereby generalizing
the expansions in \eqref{fullCMV} and \eqref{fullHess}.

\subsection{Graphical rule for the entry-wise expansion of ${\cal
S}$}

First we will present a graphical rule for predicting both the
position and the form of the non-zero entries of a snake-shaped
matrix factorization ${\cal S}$.

We will illustrate the ideas for the matrix ${\cal S}$ given by
\eqref{givenssnakeex} and Figure \ref{figsnake}. A straightforward
computation shows that the full expansion of this matrix ${\cal S}$
is given by (compare with \cite{CCG}, Example 4.5)

\begin{equation}\label{fullexpansion}
\begin{small} \hspace{-8mm}\left( \begin{array}{cccccccr}
\overline{\alpha_0} &\rho_0 &0 &0 &0 &0 &0 &0 \cdots \\
\rho_0\overline{\alpha_1} &-\alpha_0\overline{\alpha_1} &\rho_1\overline{\alpha_2} &\rho_1\rho_2 &0 &0 &0 &0 \cdots \\
\rho_0\rho_1 &-\alpha_0\rho_1 &-\alpha_1\overline{\alpha_2} &-\alpha_1\rho_2 &0 &0 &0 &0 \cdots \\
0 &0 &\rho_2\overline{\alpha_3} &-\alpha_2\overline{\alpha_3} &\rho_3\overline{\alpha_4} &\rho_3\rho_4\overline{\alpha_5} &\rho_3\rho_4\rho_5 &0 \cdots \\
0 &0 &\rho_2\rho_3 &-\alpha_2\rho_3 &-\alpha_3\overline{\alpha_4} &-\alpha_3\rho_4\overline{\alpha_5} &-\alpha_3\rho_4\rho_5 &0 \cdots \\
0 &0 &0 &0 &\rho_4 &-\alpha_4\overline{\alpha_5} &-\alpha_4\rho_5 &0 \cdots \\
0 &0 &0 &0 &0 &\rho_5\overline{\alpha_6} &-\alpha_5\overline{\alpha_6} &\rho_6 \cdots \\
0 &0 &0 &0 &0 &\rho_5\rho_6\overline{\alpha_7} &-\alpha_5\rho_6\overline{\alpha_7} &-\alpha_6\overline{\alpha_7} \cdots \\
0 &0 &0 &0 &0 &\rho_5\rho_6\rho_7 &-\alpha_5\rho_6\rho_7 &-\alpha_6\rho_7 \cdots \\
\vdots &\vdots &\vdots &\vdots &\vdots &\vdots &\vdots &\ddots
\end{array} \right) . \end{small}\end{equation}

Now the attentive reader will notice that the zero pattern of this
matrix ${\cal S}$ has some similarity with the shape of its
underlying snake as shown in Figure \ref{figsnake}. Actually, we
claim that the $(i,j)$ entry of the matrix ${\cal S}$ can be
obtained from the following recipe (the \lq E\rq\ stands for \lq
entry-wise\rq):
\begin{enumerate}
\item[E1.] Draw the snake underlying the matrix ${\cal
S}$ (cf.\ Figure \ref{figsnake});
\item[E2.] Place a right-pointing arrow on the left of the snake at height $i$;
\item[E3.] Place a left-pointing arrow on the right of the snake at height $j$;
\item[E4.] Draw the \emph{path} on the snake induced between these two arrows;
\item[E5.] If the path moves monotonically from left to right, then the $(i,j)$ entry of ${\cal
S}$ equals a product of entries of the encountered Givens
transformations
\begin{equation}\label{schurbuildingblockbis}\til{G}_{k,k+1} =
\left[\begin{array}{cc} \overline{\alpha_k} & \rho_k \\
\rho_k & -\alpha_k\end{array}\right]\end{equation}
on the path (see Step E5' below for a specification of this rule);
\item[E6.] If the path does \emph{not} move monotonically from left to right, then the $(i,j)$ entry of $S$ equals zero.
\end{enumerate}

\begin{figure}[htbp]
\begin{center}
    \subfigure[]{\label{figentrywise0}\includegraphics[scale=0.3]{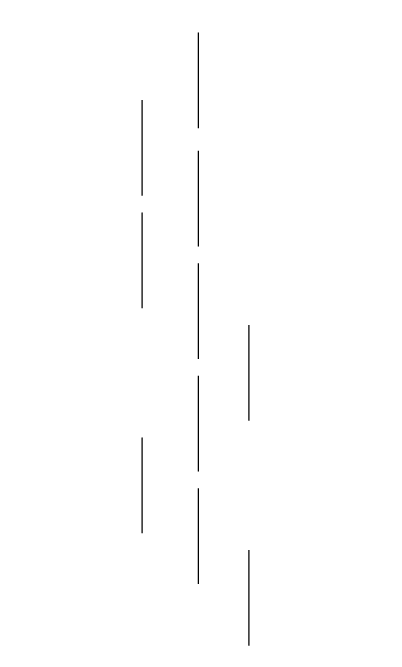}}\hspace{13mm}
    \subfigure{\label{figentrywise1brol}\includegraphics[scale=0.3]{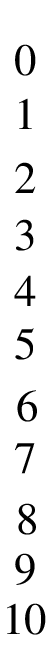}}\hspace{3mm}
    \setcounter{subfigure}{1}
    \subfigure[]{\label{figentrywise1}\includegraphics[scale=0.3]{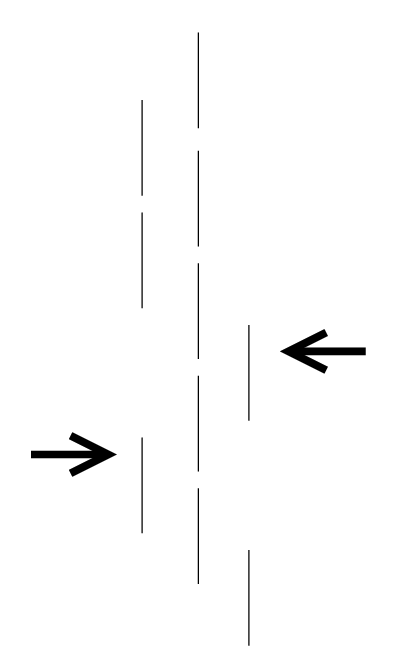}}\hspace{20mm}
    \subfigure[]{\label{figentrywise2}\includegraphics[scale=0.3]{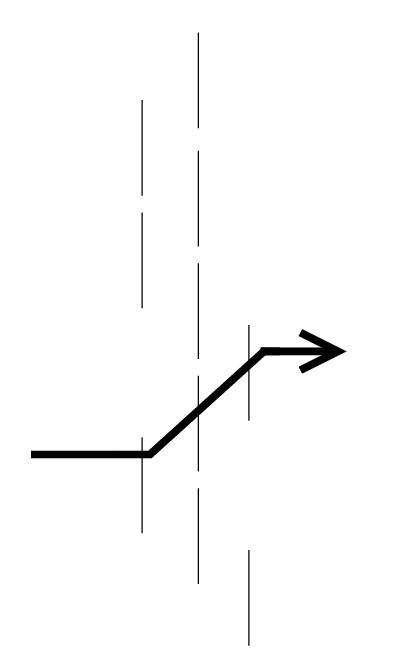}}\hspace{20mm}
\end{center}
\caption{The figure shows (a) the snake shape underlying the
matrix ${\cal S}$ in Equation \eqref{givenssnakeex}, (b) the
arrows on the left and on the right of the snake at height 7 and
5, respectively, and (c) the path on the snake induced between
these two arrows. From this information, the value of the $(7,5)$
entry of the matrix ${\cal S}$ can be determined.}
\label{figentrywise}
\end{figure}

Let us illustrate this recipe for the $(7,5)$ entry of the matrix
${\cal S}$ (recall that we label the rows and columns of this matrix
starting from the index 0). The recipe is shown for this case in
Figure \ref{figentrywise}.

Let us comment on Figure \ref{figentrywise}. Figure
\ref{figentrywise0} shows the snake shape of the matrix ${\cal S}$
(compare with Figure \ref{figsnake}), corresponding to Step E1 in
the above recipe. Figure \ref{figentrywise1} shows the arrows on the
left and on the right of the snake at height 7 and 5, respectively,
corresponding to Steps E2 and E3. The path on the snake induced
between these two arrows is shown in Figure \ref{figentrywise2},
corresponding to Step E4. Note that this path moves monotonically
from left to right and passes through the Givens transformations
$G_{7,8}$, $G_{6,7}$ and $G_{5,6}$. From Step E5 it then follows
that the $(7,5)$ entry of the matrix ${\cal S}$ is a product of
entries of these three Givens transformations. Actually, it equals
$\rho_5\rho_6\overline{\alpha_7}$ (compare with
\eqref{fullexpansion}).

\begin{figure}[htbp]
\begin{center}
    \subfigure[]{\includegraphics[scale=0.3]{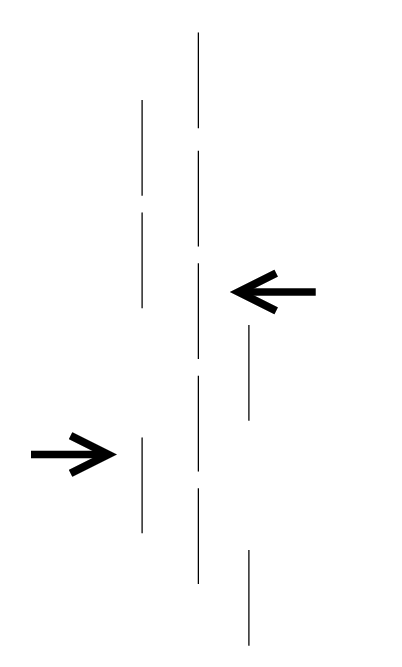}}\hspace{20mm}
    \subfigure[]{\includegraphics[scale=0.3]{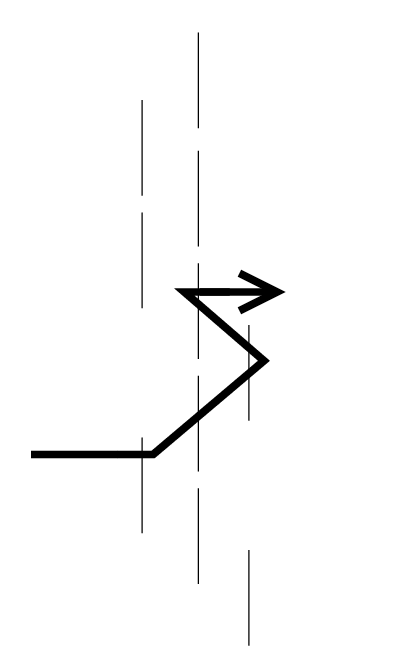}}\hspace{20mm}
\end{center}
\caption{For the matrix ${\cal S}$ in Equation \eqref{figsnake}, the
figure shows (a) the arrows on the left and on the right of the
snake at height 7 and 4, respectively, and (b) the path on the snake
induced between these two arrows. Since now the path does \emph{not}
move monotonically from left to right, it follows that the $(7,4)$
entry of ${\cal S}$ equals zero. This corresponds again with
\eqref{fullexpansion}.} \label{figentrywisebis}
\end{figure}

As a second example, let us consider the $(7,4)$ entry of the
matrix $S$. The recipe is illustrated for this case in Figure
\ref{figentrywisebis}.

Note that in the above example concerning the $(7,5)$ entry of the
matrix ${\cal S}$, we noticed that this entry equals the product
of the (complex conjugate of) the \emph{Schur parameter}
$\overline{\alpha_7}$, on the one hand, and the
\emph{complementary Schur parameters} $\rho_6$, $\rho_5$, on the
other hand. To complete our description, let us now state an
\emph{a priori} rule to determine which of the four entries in
\eqref{schurbuildingblockbis} each Givens transformation
$G_{k,k+1}$ on the path in Step E5 contributes.

Let us explain this rule for the first Givens transformation
$G_{7,8}$ through which the path in Figure \ref{figentrywise2}
passes (note that $G_{7,8}$ corresponds to the bottom leftmost line
segment on the path in Figure \ref{figentrywise2}). First, we will
determine the \emph{row index} of the entry contributed by
$G_{7,8}$. To this end, imagine that we are in the line segment
corresponding to  $G_{7,8}$ and that we move \emph{leftwards} on the
path. It is then seen from Figure \ref{figentrywise2} that we leave
this line segment through its topmost index; hence we claim that the
sought entry of $\tilde{G}_{7,8}$ will be in its topmost row.

Next, to find the \emph{column index} of the entry contributed by
$\tilde{G}_{7,8}$, imagine again that we start in the line segment
corresponding to $G_{7,8}$ but move this time \emph{rightwards} on
the path. Since the path in Figure \ref{figentrywise2} proceeds
upwards from left to right, we move then to the position of smaller
indices. Hence the sought entry of $\tilde{G}_{7,8}$ will be in its
column with the smallest index, which is column 0. We conclude that
the sought entry of $\tilde{G}_{7,8}$ lies in the $(0,0)$ position
of \eqref{schurbuildingblockbis}; this gives us
$\overline{\alpha_7}$.

The entries contributed by $G_{6,7}$ and $G_{5,6}$ can be found in a
similar way. The reader can check that in both cases, the relevant
entries of $G_{6,7}$, $G_{5,6}$ are positioned in the $(1,0)$ entry
of \eqref{schurbuildingblockbis}.

To summarize these ideas, let us introduce some notations. Denote
with $G_{r,r+1}$ and $G_{t,t+1}$ the two \emph{outermost} line
segments of the path in Step E5. Note that $r\in\{i-1,i\}$ and
$t\in\{j-1,j\}$, with the precise value of $r$ and $t$ depending
on the shape of the snake. For the example of Figure
\ref{figentrywise2}, we have $r=i=7$ and $t=j=5$.

Denote with $\mathcal{K}$ the set of indices $k$ of the
\emph{innermost} Givens transformations $G_{k,k+1}$ on the path.
Explicitly, $\mathcal{K}$ equals $\{r+1,\ldots,t-1\}$ if $r<t$ and
$\{t+1,\ldots,r-1\}$ if $r>t$ (it is understood that
$\mathcal{K}=\emptyset$ when $|r-t|=1$).

It is easily seen from the above graphical rule that the
$\{G_{k,k+1}\}_{k\in\mathcal{K}}$ in Step E5 always contribute
their \emph{complementary Schur parameter} $\rho_k$, while
$G_{r,r+1}$ and $G_{t,t+1}$ can contribute each of their entries.
In fact, we state the following specification of Step E5.

\begin{enumerate}
\item[E5'.] Under the assumptions of Step E5, and using the above
notations, the $(i,j)$ entry of the matrix ${\cal S}$ equals
\begin{equation*}
x_{r}\cdot\left(\prod_{k\in\mathcal{K}}\rho_{k}\right)\cdot y_t,
\end{equation*}
where $x_{r}\in\{ \overline{\alpha_r},\rho_r,-\alpha_r \}$ and
$y_{t}\in\{ \overline{\alpha_t},\rho_t,-\alpha_t \}$ are the
entries of $\tilde{G}_{r,r+1}$ and $\tilde{G}_{t,t+1}$ which can
be found as described in the paragraphs above\footnote{Explicitly,
$x_r$ is the $(i-r,b)$th entry of $\tilde{G}_{r,r+1}$ and $y_{t}$
is the $(1-b,j-t)$th entry of $\tilde{G}_{t,t+1}$, where the
boolean $b$ is defined by $b=0$ if $r>t$ or $b=1$ if $r<t$.}: it
suffices each time to imagine that we are in the line segment
corresponding to the current Givens transformation, and then
imagine moving leftwards or rightwards on the path, to obtain the
row and the column index in \eqref{schurbuildingblockbis},
respectively.
\end{enumerate}

\subsection{Proof of the graphical rule}

The proof that the recipe in Steps E1-E6, E5' leads to the correct
form of the $(i,j)$ entry of ${\cal S}$ follows by just expanding
the matrix ${\cal S}$ in an appropriate way. Let us sketch here the
main steps of the proof.\smallskip

\bewijs. Throughout the proof, the Givens transformations under the
$\prod$-symbol are understood to be multiplied in the order
described in Section \ref{subsectionsnakeshapes}. We will assume for
definiteness that either $i<j$ or $i=j$ and $r<t$. Consider the
given snake-shaped matrix factorization ${\cal S} =
\prod_{k=0}^{\infty} G_{k,k+1}$. Define the \lq sub-snake\rq\
\begin{equation}\label{subsnake1}{\cal S}_{i,j} := \prod_{k=i-1}^{j} G_{k,k+1}.\end{equation}

It is clear that the $(i,j)$ entry of ${\cal S}$ depends only on the
sub-snake ${\cal S}_{i,j}$. This follows since the other Givens
transformations can be considered as operations on rows and columns
$\{ 1,2,\ldots,i-1\}\cup\{ j+1,j+2,\ldots\}$ of ${\cal S}_{i,j}$;
hence indeed they cannot influence the $(i,j)$ entry of ${\cal
S}_{i,j}$.

Assume now that $s_{l}=1$ for some $l\in\{i,\ldots,j\}.$ This
means that the line segment $G_{l,l+1}$ is positioned to the left
of $G_{l-1,l}$. We can then factor \eqref{subsnake1} as
\begin{equation}\label{snaketwoparts}\left( \prod_{k=l}^{j} G_{k,k+1} \right) \cdot \left(\prod_{k=i-1}^{l-1}
G_{k,k+1}\right).
\end{equation} We distinguish between three cases:

\begin{itemize} \item Suppose that $l\in\{i+1,\ldots,j-1\}$.
The leftmost factor in \eqref{snaketwoparts} can be considered as a
row operation acting on rows $l,\ldots,j+1$ of the rightmost factor
in \eqref{snaketwoparts}. By assumption, these row indices are all
strictly larger than $i$; hence this factor cannot influence the
$(i,j)$ entry of ${\cal S}_{i,j}$. Similarly, the rightmost factor
in \eqref{snaketwoparts} acts on columns $i-1,\ldots,l$, which by
assumption are all strictly smaller than $j$.
We conclude that the $(i,j)$ entry can be influenced by
\emph{none} of the factors $G_{k,k+1}$ in \eqref{snaketwoparts},
and hence it simply equals the $(i,j)$ entry of the identity
matrix, i.e., it equals zero. This proves the conclusion in Step
E6.
\item Suppose that $l=i$. In contrast to the
previous case, we can now only conclude that the rightmost factor
$G_{i-1,i}$ in \eqref{snaketwoparts} can be removed from further
consideration. This corresponds to the fact that $r$ equals $i$
(and not $i-1$) in this case.
\item Suppose that $l=j$. Similarly as in the
previous case, we can then conclude that the leftmost factor
$G_{j,j+1}$ in \eqref{snaketwoparts} can be removed from further
consideration. This corresponds to the fact that $t$ equals $j-1$
(and not $j$) in this case.
\end{itemize}

Getting rid of all the redundant factors $G_{k,k+1}$ as described
above, we are left with either the identity matrix or with a
sequence of Givens transformations following a unitary Hessenberg
shape (cf.\ Figure \ref{fighess}). The relevant entries of this
matrix can be computed using a straightforward calculation and are
easily seen to correspond to the given rules in Steps E5 and E5'
(compare with \eqref{fullHess}). We omit further details.
  $\bol$

\subsection{Some corollaries}

A first corollary is the following.

\begin{gev} \label{bandwidthsnake}(Upper and lower bandwidth of ${\cal S}$)
The upper bandwidth of the snake-shaped matrix factorization ${\cal
S}$ equals the length of the longest sub-snake of ${\cal S}$ whose
line segments are linearly aligned in the top left-bottom right
order (cf.\ Figure \ref{fighess}). Similarly, the lower bandwidth of
${\cal S}$ equals the length of the longest sub-snake of ${\cal S}$
whose line segments are linearly aligned in the top right-bottom
left order.
\end{gev}

It follows from Corollary \ref{bandwidthsnake} that the unitary
five-diagonal matrices ${\cal C}$ have the smallest bandwidth of all
snake-shaped matrix factorizations ${\cal S}$; they have in fact
bandwidth 2 in both their lower and upper triangular part and hence
are five-diagonal.

A related result on the minimality of the matrix ${\cal C}$ is the
fact \cite{CMV3} that any infinite unitary matrix $A$ having lower
bandwidth 1 and finite upper bandwidth $n$ is \lq trivial\rq\ in
the sense that $A$ can be decomposed as a direct sum of matrices
of size at most $n+1$. This result can be shown using only some
basic linear algebra by noting that under the above conditions on
the matrix $A$, this matrix is isometric Hessenberg and hence
allows a factorization of the form \eqref{AGRfactorization}. The
condition on the upper bandwidth of $A$ then easily implies that
from each tuple of $n+1$ subsequent Givens transformations
$G_{k,k+1}$ in \eqref{AGRfactorization}, there must be at least
one for which $G_{k,k+1}$ has vanishing off-diagonal elements; we
omit further details.

A second corollary of the above results can be easily proven from
\eqref{firstbasis} and Lemma \ref{lema}. Here, the elements of the
matrix ${\cal S}$ given by \eqref{firstbasis} are expressed in terms
of the inner product \eqref{definnerproduct} and the orthonormal
Szeg\H{o} polynomials (see also Theorem 4.1 in \cite{CCG}).
\begin{gev}\label{innercoro}
By introducing the notation
$$f_i=\left\{ \begin{array}{cl}
\varphi_i(z) &\textrm{if } \;s_i=0, \\
\varphi_i^*(z) &\textrm{if } \;s_i=1,
\end{array} \right.$$
then the entries of the snake-shaped matrix factorization ${\cal
S}=\left( \eta_{i,j} \right)_{i,j\geq 0}$ are given for all $i
\geq 0$ and $k \geq 1$ by $\eta_{i,i}=\langle f_i , zf_i \rangle$
and by
$$\begin{array}{l}
\eta_{i+k,i}=\left\{ \begin{array}{cl} \langle
f_{i+k},z^{k+s_{i+k}}f_i \rangle &if \;s_{i+1}=\cdots=s_{i+k-1}=1, \\
0 &\textrm{other case},
\end{array} \right. \\
\\
\eta_{i,i+k}=\left\{ \begin{array}{cl} \langle
f_i,z^{1-s_{i+k}}f_{i+k} \rangle &if \;s_{i+1}=\cdots=s_{i+k-1}=0, \\
0 &\textrm{other case},
\end{array} \right.
\end{array}$$
where when $k=1$, the condition $s_{i+1}=\cdots=s_{i+k-1}\in\{ 0,1
\}$ is understood to be always valid.
\end{gev}
$\bol$

As a consequence of Corollary \ref{innercoro} and the graphical
rule, by choosing appropriate generating sequences one can easily
deduce a direct proof of Propositions 1.5.8, 1.5.9 and 1.5.10 in
\cite{BS}.

\section{Connection with Szeg\H{o} quadrature formulas}
\label{sectionquadrature}

In this section we describe some connections between snake-shaped
matrix factorizations and Szeg\H{o} quadrature formulas. The
results in this section are actually known for the isometric
Hessenberg and unitary five-diagonal cases, and the extension to a
general snake-shaped matrix factorization ${\cal S}$ turns out to
be rather trivial. Nevertheless, we include these results here for
completeness of the paper.

Throughout this section, we shall be dealing with a fixed measure
$\mu$ as described in Section \ref{subsectionlaurent} and we will
be concerned with the computation of integrals of the form
$$I_{\mu}(f):=\int_{\TT} f(z)d\mu(z)=\int_{-\pi}^{\pi}
f(e^{i\theta})d\mu(\theta),$$ by means of so-called {\em Szeg\H{o}
quadrature formulas}. Such rules appear as the analogue on the
unit circle of the Gaussian formulas when dealing with estimations
of integrals supported over intervals on the real line $\er$. For
a fixed positive integer $n\in\enn\setminus\{0\}$, an $n$-point
Szeg\H{o} quadrature is of the form
$$I_{n}(f):=\sum_{j=1}^{n} \lambda_j f(z_j) ,\;\; z_j \in \TT
,\;\;j=1, \ldots,n ,\;\; z_j \neq z_k \;\textrm{if} \;j \neq k,$$
where the {\em nodes} $\{ z_j \}_{j=1}^{n}$ and {\em weights} $\{
\lambda_j \}_{j=1}^{n}$ are determined in such a way that the
quadrature formulas are exact in subspaces of Laurent polynomials
whose dimension is as high as possible. The characterizing
property is that $I_n(L)=I_{\mu}(L)$ for all $L \in \textrm{span}
\{z^j : j=-n+1,\ldots,n-1\}$ (the optimal subspace): see e.g.
\cite{RO,Gragg1,Jo}, \cite[Chapter 4]{GSz}.

In what follows, we will use the notations ${\cal H}$, ${\cal C}$
and ${\cal S}$ for the isometric Hessenberg, unitary five-diagonal
and snake-shaped matrix factorization induced by the generating
sequence $\{ p_n \}_{n}$, respectively. As we have already seen,
these matrices can all be factorized as $\prod_{k=0}^{\infty}
G_{k,k+1}$, where the $G_{k,k+1}$ are canonically fixed by
\eqref{defgivens} and \eqref{schurbuildingblock}, but where the
factors under the $\prod$-symbol may occur in a certain order (cf.
Section \ref{subsectionsnakeshapes}).

We start with the following result, which seems to be
essentially\footnote{Theorem \ref{stprincipal} is not explicitly
stated in \cite{Gragg1}, but it can be easily deduced from the
results in that paper.} due to Gragg. It is the unitary analogue
of a well-known result for the Jacobi matrix when the measure
$\mu$ is supported on the real line $\er$.

\begin{st} \label{stprincipal} (Gragg \cite{Gragg1}) The eigenvalues of the principal $n \times n$
submatrix of the isometric Hessenberg matrix ${\cal H}$ are the
zeros of the $n$th Szeg\"o polynomial $\varphi_n(z)$.
\end{st}

Here with the principal $n \times n$ submatrix of ${\cal H}$ we
mean the submatrix formed by rows and columns $0$ up to $n-1$ of
${\cal H}$.

\begin{prop} \label{propprincipal} (Watkins \cite{Wat}, Cantero, Moral and Vel\'azquez \cite{Ca})
Theorem \ref{stprincipal} also holds for the unitary five-diagonal
matrix ${\cal C}$, i.e., the eigenvalues of the principal $n
\times n$ submatrix of ${\cal C}$ are the zeros of the $n$th
Szeg\"o polynomial $\varphi_n(z)$.
\end{prop}

The above results hold in fact for any snake-shaped matrix
factorization $S$; see further.

For the present discussion, a drawback of Theorem
\ref{stprincipal} and Proposition \ref{propprincipal} is that the
principal $n \times n$ submatrix of $\{ {\cal H}$, ${\cal C},
{\cal S}\}$ is in general \emph{not} unitary anymore and hence has
eigenvalues \emph{strictly} inside the unit disk. This means that
these eigenvalues are not suited as nodes for the construction of
an $n$-point Szeg\H{o} quadrature formula.

The solution to the above drawback is to slightly modify the
principal $n \times n$ submatrix of ${\cal S}$ in such a way that it
becomes unitary. Its eigenvalues will then be distinct, exactly on
the unit circle $\TT$ and turn out to be precisely the required set
of nodes.

To achieve this in practice, Gragg \cite{Gragg1} and also Watkins
\cite{Wat} introduced the idea to redefine the $(n-1)$th Givens
transformation $\tilde{G}_{n-1,n}$ by
\begin{equation}\label{decoupledgivens} \tilde{G}_{n-1,n}:=
\left[\begin{array}{cc} e^{i\theta} & 0 \\ 0 & e^{i\tilde{\theta}}
\end{array}\right],\end{equation}
where $\theta, \tilde{\theta}\in\er$ denote arbitrary parameters
(the second of them will actually be irrelevant for what follows).

With this new choice of $\tilde{G}_{n-1,n}$, we can \lq absorb\rq\
the factors $e^{i\theta}, e^{i\tilde{\theta}}$ in the previous and
next Givens transformation $G_{n-2,n-1}$ and $G_{n,n+1}$,
respectively. This means that we redefine
\begin{equation}\label{absorbsign1}
\tilde{G}_{n-2,n-1} := \tilde{G}_{n-2,n-1}\cdot
\left[\begin{array}{cc} 1 & 0 \\ 0 & e^{i\theta}
\end{array}\right],\quad \textrm{if } s_{n-1}=0,
\end{equation}
while in case $s_{n-1}=1$ we redefine $\tilde{G}_{n-2,n-1}$ by the
same formula \eqref{absorbsign1} but now with the factors multiplied
in the reverse order. Similarly, we redefine
\begin{equation}\label{absorbsign2} \tilde{G}_{n,n+1} :=
\left[\begin{array}{cc} e^{i\tilde{\theta}} & 0 \\ 0 & 1
\end{array}\right]\cdot \tilde{G}_{n,n+1},\quad \textrm{if } s_{n}=0, \end{equation}
while in case $s_{n}=1$ we redefine $\tilde{G}_{n,n+1}$ by the same
formula \eqref{absorbsign2} but now with the factors multiplied in
the reverse order. We can then put
\begin{equation} \tilde{G}_{n-1,n}:= I_2.
\end{equation}

Note that after the above updates, the value of the snake-shaped
matrix factorization ${\cal S}$ remains unchanged but we have
succeeded to transform the Givens transformation $\tilde{G}_{n-1,n}$
in \eqref{decoupledgivens} into the identity matrix $I$. Then it is
easily seen that the snake shape of ${\cal S}$ can be \lq broken\rq\
into two pieces, in the sense that ${\cal S} = UV$ where $U =
\prod_{k=0}^{n-2} G_{k,k+1}$ is the submatrix formed by rows and
columns $0,\ldots,n-1$ of ${\cal S}$, while $V =
\prod_{k=n}^{\infty} G_{k,k+1}$ is the submatrix formed by rows and
columns $n,\ldots,\infty$. Note that the matrices $U$ and $V$ have a
complementary zero pattern and hence they commute with each other.

Let us now denote with ${\cal S}_{n-1}:=\prod_{k=0}^{n-2}
G_{k,k+1}$ the topmost part of the \lq broken\rq\ snake ${\cal
S}$. Note that ${\cal S}_{n-1}$ is a snake-shaped matrix
factorization of size $n\times n$; in particular it is still
\emph{unitary}. Note also that this matrix depends on the
parameter $\theta\in\er$ by means of \eqref{absorbsign1}.

\begin{opm}\label{remarkprincipal} The principal $n\times n$ submatrix of ${\cal S}$ can be obtained in the same
way as above, but now replacing the role of $e^{i\theta}\in\TT$ in
\eqref{decoupledgivens} by the original matrix entry
$\overline{\alpha_n}$. Note however that $\overline{\alpha_n}$ lies
strictly inside the unit disk and hence the resulting $n\times n$
submatrix is \emph{not} unitary anymore; cf.\ the motivation earlier
in this section.
\end{opm}

One has then the following result.

\begin{st} \label{stquad} (Gragg \cite{Gragg1}) Let $\theta\in\er$ be fixed.
Using the above construction, the eigenvalues
of ${\cal H}_{n-1}$ are distinct, belong to $\TT$ and appear as
nodes in an $n$-point Szeg\H{o} quadrature formula for the measure
$\mu$. The corresponding quadrature weights are the first
components of the normalized eigenvectors of ${\cal H}_{n-1}$.
\end{st}

\begin{prop} (Watkins \cite{Wat})
Theorem \ref{stquad} also holds for the matrix ${\cal C}_{n-1}$.
\end{prop}


Here with \lq normalized\rq\ eigenvectors we mean that the
eigenvectors should be scaled in such a way that they form an
orthonormal system and that their first components are real
positive numbers.

The characteristic polynomial of the above matrix ${\cal H}_{n-1}$
(or equivalently, ${\cal C}_{n-1}$) is known as a monic
\emph{para-orthogonal polynomial} of degree $n$ \cite{Jo}. Note that
this polynomial depends on the free parameter $\theta$, and hence
there is in fact a one-parameter \emph{family} of para-orthogonal
polynomials (and so, a one-parameter family of $n$-point Szeg\H{o}
quadrature formulas for $\mu$).

Now one could ask why there is such a similarity between ${\cal
H}$ and ${\cal C}$ in the above results. This is explained by the
following basic observation, which is essentially due to Ammar,
Gragg and Reichel \cite{AGR} for the case of ${\cal H}$ and ${\cal
C}$.

\begin{prop}\label{propAGR}(Based on Ammar, Gragg and Reichel \cite{AGR}) Let $\theta\in\er$ be fixed. Then the eigenvalues and
the first components of the normalized eigenvectors of ${\cal
S}_{n-1}$ depend on the Schur parameters but \emph{not} on the shape
of the snake.
\end{prop}

\bewijs. Recall that the snake-shaped matrix factorization is given
by ${\cal S}_{n-1} = \prod_{k=1}^{n-2} G_{k,k+1}$, for some order of
the factors. But it is a general fact that the matrices $AB$ and
$BA$ have the same eigenvalues; this follows from the similarity
transformation \begin{equation}\label{similarity} AB\mapsto
A^{-1}(AB)A = BA.\end{equation} By applying this idea recursively
for the choice $A=\prod_{k=l}^{n-2} G_{k,k+1}$, for $l=n-2,\ldots,1$
(only those indices $l$ for which $s_l=1$ have to be treated), one
can succeed to rearrange the Givens transformations of $S_{n-1}$
into the unitary Hessenberg form $G_{0,1}G_{1,2}\cdots G_{n-2,n-1}$
(compare with \eqref{AGRfactorizationbis}). It follows that the
eigenvalues of ${\cal S}_{n-1}$ are indeed independent of the order
of the factors $G_{k,k+1}$, i.e., they are independent of the shape
of the snake.

The same argument also shows that the first components of the
normalized eigenvectors are independent of the shape of the snake.
To see this, consider the eigen-decomposition ${\cal S}_{n-1} =
UDU^{*}$, where $D$ is a diagonal matrix containing the eigenvalues,
and $U$ is a unitary matrix whose columns are the eigenvectors,
scaled in such a way that the first row of $U$ has real positive
entries. The point is now that the only Givens transformation of
${\cal S}_{n-1}$ acting on the $0$th index is $G_{0,1}$; but in the
above argument the latter can only appear as the $B$-factor in
\eqref{similarity}, and hence the first row of $U$ is easily seen to
remain unchanged under the similarity \eqref{similarity}. $\bol$

\begin{gev} Theorems \ref{stprincipal} and \ref{stquad} hold with ${\cal H}$ replaced by
any snake-shaped matrix factorization ${\cal S}$.
\end{gev}

\bewijs. This follows from Proposition \ref{propAGR} and Remark
\ref{remarkprincipal}.$\bol$

Note that Proposition \ref{propAGR} implies that the eigenvalue
problems for the matrices ${\cal H}_{n-1}$, ${\cal C}_{n-1}$ and
${\cal S}_{n-1}$ are conceptually equivalent. Interestingly, these
problems turn out to be also \emph{numerically} equivalent since,
for reasons of efficiency and numerical stability, the eigenvalue
computation for $\{{\cal H}_{n-1}$, ${\cal C}_{n-1}$, ${\cal
S}_{n-1}\}$ should preferably be performed using their factorization
as a product of Givens transformations, rather than using their
entry-wise expansions.

Finally, we mention that the development of extensions of
Szeg\H{o} quadrature formulas and the investigation of the
connection between them and Gauss quadrature formulas on the
interval $[-1,1]$ are active areas of research: see e.g.
\cite{Bul,Cr,RO,JaRe} and references therein found. A whole
variety of practical eigenvalue computation algorithms for unitary
Hessenberg and five-diagonal matrices has already been developed
in the literature. In \cite{Rut}, Rutishauser designed an
LR-iteration. Implicit QR-algorithms for unitary Hessenberg
matrices were described and analyzed in
\cite{q232,q352,m281,m240}. In \cite{ARS,GGCC} and the references
therein, divide and conquer algorithms were constructed. Other
approaches are an algorithm using two half-size singular value
decompositions \cite{AGR}, a method involving matrix pencils
\cite{BGE}, and a unitary equivalent of the Sturm sequence method
\cite{BGH}.

\section*{Acknowledgment} The authors thank professors A.~Bultheel and
A.~Kuijlaars for useful suggestions and the referees for valuable
comments.

\end{document}